\documentclass[10pt]{article}
\usepackage{times}
\usepackage{amsmath}
\usepackage{amssymb}
\usepackage{amscd}
\usepackage{stmaryrd}
\usepackage{xspace}
\input xy
\xyoption{all}
\usepackage[vcentermath]{youngtab}

\newtheorem{thm}{THEOREM}
\newtheorem{lemma}{LEMMA}
\newtheorem{prop}{PROPOSITION}
\newtheorem{cor}{COROLLARY}
\newtheorem{deff}{DEFINITION}
\newcommand{\beq}{\begin{equation}}
\newcommand{\eeq}{\end{equation}}
\newcommand{\ba}{\begin{array}}
\newcommand{\ea}{\end{array}}
\newcommand{\beqa}{\begin{eqnarray}}
\newcommand{\eeqa}{\end{eqnarray}}

\def\S{{\mathfrak S}}

\def\hn{{\mathcal H}_n(q)}
\def\hr{{\mathcal H}_r(q)}

\def\h{{\mathcal H}(q)}

\def\Pl{{\mathfrak{P}}}
\def\PP{{\mathfrak{PP}}}
\def\prPP{{\PP^{Lie}}}

\def\Uq{{U_q{\mathfrak{gl}}}}
\def\C{{\mathbb C}}

\def\Z{{\mathbb Z}}
\def\Q{{\mathbb Q}}
\def\S{{\mathfrak S}}
\def\I{{\mathfrak I}}

\begin{document}

\title{Pre-Plactic Algebra and Snakes}

\author{Todor Popov \\
\it Institute for Nuclear Research and Nuclear Energy\\ 
\it Bulgarian Academy of Sciences \\ 
\it Sofia, BG-1784, Bulgaria \\
 tpopov@inrne.bas.bg
}

\maketitle

\abstract{

We study a factor Hopf  algebra $\PP$ 
 of the Malvenuto-Reutenauer  convolution algebra 
of functions on symmetric groups $\S=\oplus_{n\geq 0}\C[\S_n]$
that we coined {\it pre-plactic algebra}. 
The pre-plactic algebra  admits  the  Poirier-Reutenauer algebra  based on Standard Young Tableaux 
 as a factor
 and it is  closely related to the quantum pseudo-plactic algebra  introduced by Krob and Thibon in the non-commutative character theory of  quantum group comodules.
The connection between the
quantum pseudo-plactic algebra  and the pre-plactic algebra   is similar to
the connection between the Lascoux-Sch\"utzenberger plactic algebra 
 and the  Poirier-Reutenauer  algebra.
We show that the  dimensions of the pre-plactic algebra are given by the numbers of alternating permutations (coined snakes after V.I. Arnold). 
Pre-plactic algebra is instrumental in calculating the Hilbert-Poincar\'e series of the quantum pseudo-plactic algebra.}

\section{Quantum Pseudo-Plactic Algebra}

The diagonal matrix  elements of the coorinate ring
$\C[GL_q(V)]$ of a quantum group close a subalgebra which is the noncommutative reincarnation of the  algebra of the functions on the torus. Krob and Thibon conjectured \cite{KT}  that the diagonal subalgebra 
in $\C[GL_q(V)]$ is a cubic algebra that will be referred to as the {\it quantum pseudo-plactic algebra}.

\begin{deff} 
The quantum  pseudo-plactic algebra $\PP_q(W)$  
of the vector space $V$ with ordered basis,  
$W = \bigoplus_{i=1}^D \C(q) a_i$ is the 
quotient of the 
free associative algebra $T(W)$ of $W$
by the two-sided ideal $\I$ generated
 by the cubic relations
\beq
\label{ppq}
\ba{lccc}
[[a,c],b] &=&0  &a<b<c \\[4pt]
[[a,b],a]_{q^2} &=&0  &a<b \\[4pt]
[b,[a,b]]_{q^2} &=&0  &a<b 
\ea
\eeq
which we are referring to as quantum pseudo-Knuth relations.
Here $[x,y]_q$ stands for the deformed commutator
$[x,y]_q:=xy-q yx$ and $[x,y]=[x,y]_1$. Let $\I_q$ be
the  two-sided ideal  generated by the
quantum pseudo-Knuth 
relations then one has
\[
\PP_q(W) = T(W)/\I_q(W) \ .
\]
\end{deff}
{\bf Remark.}	The algebra		$\PP(W)$ has as a factor the coordinate ring of the quantum torus,  $y x=q^2 x y$ for $x<y$ .

The specialization at $q=1$ of the quantum pseudo-plactic algebra
relations 
\beq
\label{}
\ba{cccc}
[[a,c],b] &=&0  &a<b<c \\[4pt]
[[a,b],a] &=&0  &a<b \\[4pt]
[b,[a,b]] &=&0  &a<b 
\ea
\eeq
 involve only commutators  thus $\PP_1(W) = T(W)/\I_1(W)$ is an universal enveloping algebra of a Lie algebra.
The specialization $\PP_1(W)$ will be referred to as pseudo-plactic algebra and denoted by $\PP(W) =T(W)/\I(W)$.

\noindent
{\bf Remark.} The specialization at $q=1$ of the  diagonal subalgebra  $\C[GL_q(V)]$  is a commutative algebra which is different from  the cubic pseudo-plactic algebra $\PP(W)$.


The quantum Schur-Weyl duality  is 
the  double centralizing property  on ${\rm{End}}(T(V))$ of the action of the Hecke algebra $\h=\bigoplus_{r \geq 0} \hr$ and the (co)action
of the quantum group $\C[GL_q(V)]$. A
Schur functor is mapping a $\h$-module to a $\C[GL_q(V)]$-comodule
(or equivalently to a $\Uq(V)$-module).
The quantum pseudo-Knuth ideal $\I_q(W)$ inherts 
 $\C[GL_q(V)]$-comodule structure upon restriction from 
the diagonal subalgebra of $\C[GL_q(V)]$. Thus
 the quantum pseudo-plactic algebra $\PP_q(W)$ 
allow for a description \cite{P} based on Schur functors
\[
\PP_q(W)=\bigoplus_{n\geq 0} \PP_q (W)_n = 
\bigoplus_{n\geq 0} \PP_q (n)\otimes_{\hn} W^{\otimes n}  
\]
where $W$ is the diagonal in  
$V\otimes V^\ast$, {\it i.e.}, the elements
invariant under transposition. 
Let  $\h$-module $\I_q=\oplus_{r\geq 0}\I_q(r)$ be the Schur functor 
pre-image of the quantum pseudo-Knuth ideal $\I_q(W)$.
Then the
the collection $\{\PP_q(r) \}_{r\geq 0}$ of  $\hr$-modules determines a $\h$-module
\[
\PP_q=\bigoplus_{r\geq 0} \PP_q(r) = \bigoplus_{r\geq 0} \hr/\I_q(r)  \ 
\]
which will be referred to as quantum pre-plactic algebra.

Similarly  we define the pre-plactic algebra $\PP$ to be a collection $\{ \PP(r) \}_{r\geq 0}$ of  $\S_r$-modules (see Definition \ref{prep}).
We are going to endow the pre-plactic algebra $\PP$   with a structure of a Hopf algebra 
induced from the Malvenuto-Reutenauer 
Hopf structure on 
$\S=\bigoplus_{n\geq 0} \Q\S_n $ which we now briefly revise.

\section{Malvenuto-Reutenauer Hopf Algebra}

Malvenuto and Reutenauer \cite{MR}
have considered a Hopf algebra structure on $\S$ 
related to the Solomon descent
algebra \cite{S}.
There is a dual Hopf structure on $\S$ 
known also as
the algebra of free quasi-symmetric functions \cite{DHT}.

Given two permutations $\alpha \in \S_{r}$ and $\beta \in \S_{p}$, the product $\ast$ on $\S=\cup_{n\geq 0} \Q[\S_n]$ is induced (by linearity) from
$$
\begin{array}{rccc}
\ast :& \Q[\S_r]\times \Q[ \S_p] &\rightarrow & \Q[ \S_{r+p}]
\uparrow^{\S_{r+p} }_{\S_r\times \S_p}\ , \\[6pt]
&\alpha \times \beta &\mapsto&\alpha \ast \beta = \mathop{\sum}\limits_{\tiny\ba{c}st(u)=\alpha \\  st(v)=\beta \ea} 
u v  \ .
\end{array}
$$
Here $st(w)$  is the standard word  in $ \S_n$ of
the word $w=w_1 \dots w_n$ of length $n$, i.e., the image of the unique increasing injective function from 
$\{w_1, \dots w_n \}$ to $\{1, \dots  n \}$.

Similarly the coproduct $\Delta$ in  $\S$ is induced from
$$
\begin{array}{rccc}
\Delta:& \Q[ \S_r ] &\longrightarrow & 
\mathop{\bigoplus}\limits_{i=0}^{r}
\Q[\S_i]\otimes \Q[\S_{r-i}]\downarrow^{\S_r}_{\S_i\times \S_{r-i}} \ ,\\
[4pt]
&\alpha &\mapsto& \Delta(\alpha) = \sum_{i=0}^{r} \alpha_{\{1,\ldots,i\}}\otimes
st(\alpha_{\{i+1,\ldots,r\}})
\end{array}
$$
Here $\alpha_I$ is the word obtained from $\alpha$ by erasing the letters outside $I$.


Let $\langle \,\, | \, \, \rangle$ stands for  the scalar product   such that
 the  basis of permutation words in $\S$ is orthonormal $\langle \alpha | \beta \rangle=\delta_{\alpha,\beta}$, $\alpha \in \S_{r}$ and $\beta \in \S_{p}$.
The dual Hopf structure  $\ast', \Delta'$
on $\S$
 is in linear duality with the structure $\ast, \Delta$,
{\it i.e.}, it is a convolution algebra  defined through
$$
\langle \alpha \ast \beta | \gamma \rangle =
\langle \alpha \otimes \beta| \Delta' \gamma 
\rangle 
\qquad
\langle \alpha \ast' \beta| \gamma \rangle =
\langle \alpha \otimes \beta| \Delta \gamma \rangle
 \qquad \alpha, \beta, \gamma \in \S \ .
$$
giving rise to the product and the coproduct
\beqa
\alpha \ast' \beta &=& sh(\alpha, \bar{\beta}) \ , \\
\Delta'\alpha&=&\mathop{\sum}\limits_{\alpha=uv} st(u)\otimes st(v) \ .
\eeqa
The dual Hopf algebra $(\S,\ast^\prime, \Delta^\prime)$
is isomorphic to  the algebra of the free quasi-symmetric functions $\bf FQSym$ \cite{DHT}.
Throughout this paper we will be dealing only  with the bialgebra
structures of $\S$, for the explicit formulas for the antipode see
\cite{AS}.


\section{Poirier-Reutenauer Hopf algebra}
The celebrated Robinson-Schensted algorithm \cite{K} establishes a one-to-one correspondence between  permutations $\alpha$ and pairs $(P,Q)$ of Standard Young Tableaux (SYT) of the same form
$$
RS: \S \rightarrow SYT\times SYT \qquad  RS: \alpha \mapsto (P(\alpha),Q(\alpha)) \ .
$$
Two  permutations $\alpha$ and $\beta$ are $P$-equivalent, i.e., $P(\alpha)=P(\beta)$,
if and only if their words are congruent $\alpha \sim \beta$ with respect to the plactic congruence generated by the Knuth equivalence relations \cite{K} (without repeating letters)
\beq
bac \sim bca \qquad  acb\sim cab  \qquad a<b<c \ .
\eeq

\begin{deff}Let
$Tab$ be  the quotient $ Tab:= \S / \mathfrak K$
where
 $\mathfrak K$ is the ``Knuth'' submodule of $\S$ generated by the elements $\alpha -\beta$,
$\alpha,\beta \in \S$, such that $\alpha \sim \beta$.
\end{deff}
To any $t\in SYT$ we can associate its plactic class $c(t)$, which is the sum of
the elements in $\alpha\in \S$, such that $P(\alpha)=t$
\beq
\label{plcl}
c(t) = \sum_{\alpha: t=P(\alpha)} \alpha \,\,  \subset \S \ .
\eeq
\begin{deff}
Let us denote by $Tab' = \bigoplus_{t\in SYT} \, \Q \, c(t)$ the $\Q$-span of all plactic classes.
\end{deff}

Poirier and Reutenauer proved that the Malvenuto-Reutenauer Hopf structures on $\S$
induce Hopf structures on the $\S$-modules $Tab$ and $Tab'$.

\begin{thm}\cite{PR}.
The submodule $\mathfrak K$ of $\S$ is an ideal and a coideal of $(\S,\ast,\Delta)$.
The quotient $ Tab $ 
is a  Hopf\, factor-algebra $(Tab, \ast, \Delta)$.
\label{tabiso}
\end{thm}

\begin{thm}\cite{PR}. The submodule $Tab' \subset \S$ is a Hopf subalgebra $ (Tab',\ast',\Delta')$
of the Hopf algebra $ (\S,\ast',\Delta')$.
The Hopf algebras $(Tab,\ast,\Delta)$  and $ (Tab',\ast',\Delta')$ are canonically dual.
\end{thm}

The Hopf subalgebra $ (Tab',\ast',\Delta')$ of the Malvenuto-Reutenauer algebra $(\S,\ast',\Delta')$ is isomorphic to  the algebra
of the free symmetric functions $\bf FSym$ \cite{DHT}.

\section{Pre-Plactic Hopf Algebra}

Let $u\in \S_p$ and $v\in \S_r$ be (the words of) two permutations.
We denote by $\mathfrak I$ the  two-sided ideal $\I$   in $\S$ generated by elements 
\beq
\label{}
u[b[ac]]v:=u( bac - bca - acb + cab)v \in \Z[\S_{p +3+r} ] 
\qquad a<b<c
\eeq
The relations in the ideal are among the
relations of the quantum pseudo-plactic algebra $\PP_q(W)$,
namely the multilinear part of the ideal $\I_q(W)$ and $\I(W)$.
The idea is to have a kind of  ``operadic approach" to the 
algebra $\PP_q(W)$.

\begin{lemma} The submodule $\mathfrak I\subset \S$ is a Hopf ideal of the Hopf algebra $(\S, \ast, \Delta)$.
\end{lemma}
{\bf Proof.}
The ideal $\mathfrak I$ is generated in $\Z[\S_3] \subset \S$
by the element $[2[13]]$. By definition it is a graded ideal
$\mathfrak I=\oplus_{n\geq 0} \mathfrak I(n)$ with degrees
\beq
\label{ideal}
\mathfrak I(n)= \sum_{i,j:i+j+3=n} \Q[ \S_i] \ast [2[13]] \ast \Q[  \S_j] \subset \Q[\S_n] \qquad \Rightarrow \qquad
\mathfrak I= ([2[13]]) \ .
\eeq
 The module $\mathfrak I$ is a submodule of the  module $\mathfrak K$, having as generators differences of the ``Knuth" generators in
$\mathfrak K$ hence the inclusions of submodules $\mathfrak I \subset \mathfrak K \subset \S$. It is easy to see that 
$\mathfrak I$ is also coideal with respect to $\Delta.$ 
{${\Box}$}

\begin{deff}
\label{prep}
 The pre-plactic algebra $\PP$ is
the quotient Hopf algebra $\S/\mathfrak I$ 
$$\PP:=\S/ ([2[13]]) \ .$$
\end{deff}

The factor algebra $\PP$ inherits the grading of $\S$,
\[
\PP = \bigoplus_{n\geq 0} \PP(n) = \bigoplus_{n\geq 0} \Q[\S_n ]/ \mathfrak I(n) \ .
\]

We are going to calculate the dimensions $\dim \PP(n)$.
To this end we  take up with the dimensions of the graded ideal
$\mathfrak I$. We first introduce some notations about partitions.

Let $\lambda=(\lambda_1, \ldots, \lambda_k)$ be a 
partition $\lambda_1\geq \lambda_2\geq \ldots \geq\lambda_k>0$ with
 $r$ nonzero parts.  The graphic representation of $\lambda$ is the Young diagram
whose $i$-th row is of length $\lambda_i$, $i=1,\ldots,k$.
By $\lambda'$ we denote the diagram transposed to
$\lambda$. 
A dual description $\lambda=[1^{m_1}2^{m_2}\ldots 
k^{m_k}]$ is provided by the multiplicities $m_i(\lambda)=\lambda'_i-\lambda'_{i+1}$  yielding the number of rows of length $i$. One has $|\lambda|=\sum_{i=1}^k i \, m_i(\lambda)=
\sum_{i=1}^k \lambda_i$ with $k=\sum m_i(\lambda)$.

\begin{prop} The dimensions of the  graded ideal $\mathfrak I=\bigoplus_{n\geq 0} \mathfrak I(n)$ are counted by
\beq
\label{my}
\dim \mathfrak I(n)=
{\sum_{ k >0 }} \,\,
 {\sum_{\scriptsize\ba{c}{\Lambda}\vdash
 n \\  \Lambda \neq 1^n\ea}} (-1)^{1+(|{\Lambda}|-k)/2}
\frac{k!}{ \prod_{i=1}^k{m_i({\Lambda}) !}} 
\frac{n!}{\prod_{i=1}^k  {\Lambda}_i!} 
\eeq
where the sum runs 
over all partitions ${\Lambda}=({\Lambda}_1, \ldots, {\Lambda}_k)$ of $n$ with $k$ odd parts ($\Lambda_i=2\lambda_i+1$)  
except the partition $\Lambda=1^n$.

\end{prop}
{\bf Proof.}
Our strategy in counting the dimension
$\dim \mathfrak I(n)$ is a recurrent use of the formula for the dimension of an intersection  
\beqa
\label{inter}
\dim \sum_{\alpha} A_\alpha&=& \sum_{\alpha} \dim A_\alpha - \sum_{\alpha_1<\alpha_2}
\dim A_{\alpha_1}\cap A_{\alpha_2}  + \sum_{\alpha_1<\alpha_2<\alpha_3}
\dim A_{\alpha_1}\cap A_{\alpha_2} \cap A_{\alpha_3}  \nonumber
\\ &-&  \ldots
+ (-1)^{l+1} \dim A_{\alpha_1}\cap A_{\alpha_2}\cap A_{\alpha_3} \cap \ldots
 \cap A_{\alpha_l} 
\eeqa
where the  family of spaces $\{ A_\alpha\}$ 
such that 
$\mathfrak I(n) = \sum_{\alpha} A_\alpha$
is indexed by a finite ordered set
of cardinality  $l$.

We put  into correspondence  the composition $c=1^i31^j$ 
of $n=i+j+3$
with the subspace $A_c=\Q[ \S_i] \ast [2[13]] \ast \Q[  \S_j]\subset \mathfrak I(n)$. 
It is convenient  to set $A_1:=\Q[\S_1]$ then one has
$$A_{1^j}
= A_1^{\ast j}=  \Q[\S_j] \qquad A_3=\Q [2[13]]
$$
where $A_3$ is identified with the generator of the ideal $\mathfrak I= ([2[13]])$. 
With the new conventions the   expression (\ref{ideal}) for the
graded  ideal takes the succinct form
\beq
\label{ideal1}
\mathfrak I(n)=
\sum_{i=0}^{n-3} A_{1^j} \ast A_3 \ast A_{1^{n-3-i}}=\sum_{i=0}^{n-3}  A_{1^i31^{n-3-i} }  \ .
\eeq

\begin{lemma}
\label{zeroint}
The neighboring subspaces of $\mathfrak I(n)$
in the sum eq.(\ref{ideal1}) have zero intersection
\beq
\label{gint}
A_{1^k31^{l+1}}\cap A_{1^{k+1}31^l}=\emptyset \ .
\eeq
\end{lemma}
{\bf Proof of the lemma.}
The  $\S_4$-module $\mathfrak I(4)$ is a sum of two subspaces $\mathfrak I(4)=A_{31}+A_{13}$ and its dimension
is given by the intersection formula eq. (\ref{inter}) 
\[
\dim \mathfrak I(4)=\dim A_{31}+ \dim A_{13} -
\dim  A_{13}\cap A_{31}  \ .
\]
One can write explicitly the basis of $A_{31}=A_3 \ast A_1$,
it consists of the four elements in $\Q [\S_4]$
\[
r_1 =[[13]2]4 \qquad r_2=[[14]2]3 \qquad r_3=[[14]3]2 \qquad 
r_4=[[24]3]1 
\]
whereas the basis of $A_{13}=A_1\ast A_3$ is spanned by 
\[
l_1=4[[13]2] \qquad l_2=3[[14]2] \qquad l_3=2[[14]3] \qquad l_4=1[[24]3] \ .
\]

An element $w$ belongs to the  intersection
$A_{13} \cap A_{31}$ when it can be represented  
$w= \sum_i c_i l_i=\sum_i d_i r_i $
with  coefficients $c_i$ and $d_i$. 
However the only solution turns out to be $c_i=0=d_i$.
We conclude that the intersection is
empty, 
\[A_{13}\cap A_{31}=\emptyset \ . \]

The $\S_5$-module $\mathfrak I(5)$ is  in the linear
enveloppe of three subspaces,
$$\mathfrak I(5)=A_{1^23}+A_{131}+A_{31^2}$$
and the vanishing $A_{13}\cap A_{31}=\emptyset$
readily implies 
 the zero intersections $A_{1^23}\cap A_{131}=\emptyset$
and $A_{131}\cap A_{31^2}=\emptyset$.

The statement eq.(\ref{gint}) follows by induction of the degree $\mathfrak I(n)$
of the ideal $\mathfrak I$. $\Box$

However the intersection $A_{1^23}\cap A_{31^2}$ is not trivial
(see Lemma \ref{int} below).

Let us  introduce some concise notations for the intersections
$$
A_{1^23}\cap A_{31^2} =:A_5\qquad
A_{31^3}\cap A_{1^231} =: A_{51} \qquad
A_{31^3}\cap A_{1^33} =: A_{33} \ .  $$
Then the latter lemma implies the succinct expressions of
the dimensions
\beq
\ba{ccl}
\dim \I(i) &= &0 \qquad i<3 \\
\dim \I(3) &= &\dim A_3 \\
\dim \I(4) &= &\dim A_{31}+ \dim A_{13}\\
\dim \I(5)&=&\dim A_{31^2}+\dim A_{131} +\dim A_{1^23}
- \dim A_5 \\
\dim \I(6)&=&
\dim A_{31^3}+ \dim A_{131^2} + \dim A_{1^231} + \dim A_{1^33}
 \\ &&
-\dim A_{51} - \dim A_{15} - \dim A_{33} 
 \ea
\eeq
A non-trivial (maximal) intersection
$A_7:= A_{31^4} \cap A_{1^231^2} \cap A_{1^43}$
will appear in $\I(7)$. 
In general on any odd degree $2k+1$  an intersection
of $k$ subspaces
takes place.

\begin{lemma}
\label{int}
Let the space   $A_{2n+1}$ be the  intersection of $n$ 
subspaces\footnote{
Note that only next to neighboring terms from the sum
$\I(2n+1)=\sum A_{1^k3 1^{2n-2-k}}$
appear in the intersection, as the intersections
of neighboring terms is empty due to lemma \ref{zeroint}.}
$A_{1^k31^{2n-k-2}}$
\beq
\label{odd}
A_{2n+1}:=A_{31^{2n-2}}\cap A_{1^231^{2n-4}}\cap
A_{1^431^{2n-6}}\cap \ldots \cap A_{1^{2n-2}3} \in \I(2n+1) \ .
\eeq
The intersection $A_{2n+1}$ is maximal in $\I(2n+1)$
and it is one dimensional
\[
\qquad \dim A_{2n+1}=1 \ .
\] 
\end{lemma}
{\bf Proof of the lemma.} For $n=1$ the 
spatement of the lemma is trivial 
\[
\dim \I(3) = \dim A_3 =1 \ .
\]

The first non-trivial statement
is for the space $A_5$.
For a commutative algebra with three generators $x,y,z$
 one has the antisymmetric combination
\beq
\label{xyz}
x\wedge y\wedge z:=x[yz]+y[zx]+z[xy]=[yz]x+[zx]y+[xy]z
\eeq
which spans a one dimensional space. Note that
the generators $x=[15]$, $y=[24]$ and $z=3$ in $\PP(5)$
commute\footnote{ The relation
$[[15],[24]] = [[[15] 2]4]+ [2 [[15]4]]$ is implied by the Jacoby identity.} and  generate a maximal commutative subalgebra.
Hence we have a one dimensional subspace
\[
[15]\wedge[24]\wedge 3 \in \I(5) \ .
\]


The following element belongs to the intersection $A_5=A_{1^23}\cap A_{31^2}$
\[
[15][[24]3]+ [24][3[15]] \equiv 
[[24]3][15]+[3[15]][24] =[15]\wedge[24]\wedge 3
\quad mod \quad [[15][24]] \ .
\]
The relation 
$[[15],[24]] $ can be represented also as $[15]\wedge[24]$
thus for the intersection space we get
\[
A_5=A_{1^23}\cap A_{31^2}= \Q [15]\wedge[24]\wedge 3
\quad mod \quad [15]\wedge[24]  \ .
\]
Therefore $A_5$ belongs to the ideal $\mathfrak I(5)$ and  it is one dimensional
\[
\dim A_5=1 \ .
\]
So we have proven the lemma when $n=2$.

We proceed by induction.
 In $\mathfrak I(6)$  the maximal commutative subalgebra
is generated by $[16]$, $[25]$ and $[34]$. Hence
according to eq.  (\ref{xyz}) we get the antisymmetrized combination
$[16]\wedge [25]\wedge [34] \in \mathfrak I(6)$
which generates between others the following relations in $\PP(7)$,
\[
4[17]\wedge [26]\wedge [35] \in \I(7) \qquad [17]\wedge [26]\wedge [35]4 \in \I(7) \ .
\]
The maximal commutative subalgebra is $\PP(7)$ is generated by
the ``long'' generators
$[17]$, $[26]$, $[35]$ and the ``short'' generator  $4$. 
The commutative subalgebra
gives rise to  the space  
\[
[17]\wedge[26]\wedge[35]\wedge 4 \quad mod \quad 
[17]\wedge[26]\wedge[35]
\]
which belongs to 
the one-dimensional  intersection $A_7=A_{31^4}\cap A_{1^231^2}\cap A_{1^43}$.

By induction on $n$  we obtain an unique element 
in the $\S_{2n+1}$-module $\mathfrak I(2n+1)$
\[
[1 \,2n+1]\wedge [2 \, 2n] \wedge \ldots \wedge [n\, n+2]
\wedge 
n+1    \ ,
\]
which after factoring gives rise  to
the intersection space (\ref{odd})
\[
A_{2n+1}=\Q [1 \,2n+1]\wedge [2 \, 2n] \wedge \ldots \wedge [n\, n+2]
\wedge 
n+1 \quad mod \quad
[1 \,2n+1]\wedge [2 \, 2n] \wedge \ldots \wedge [n\, n+2]
\ .
\]
We are done with the lemma. $\Box$

From the above examples the general pattern for $\dim \mathfrak I(n)$ is clear, we have summation with alternating signs  over
the dimensions of subspaces $A_c =
A_{ c_1}\ast \ldots \ast A_{ c_k}\subset \Q[ \S_n]$
 labelled by  compositions $c$ of $n$ having only odd parts 
\[
\dim \mathfrak I(n) = \sum_{c: |c|=n} \pm \dim A_c
\]
where the sum is over all compositions of $n$ with odd parts
(except $c=1^n$) and the sign $\pm$ is $(-1)^{k+1}$ if
$A_c$ stems from intersection of $k$ elementary subspaces $A_{1^i31^{n-3-i}}$.

Subspaces $A_{c}$ and $A_{c'}$ whose compositions  after ordering are mapped to the same partition $\Lambda$
will have same dimension $\dim A_{c} =\dim  A_{c'}= \dim A_{\Lambda}$. Therefore the sum over compositions $c$ can be replaced
by the sum over partitions $\Lambda$ at the expense of the multiplicities
\[
\dim \mathfrak I(n) = \sum_{k> 0}
{\sum_{\scriptsize\ba{c}{\Lambda}\vdash
 n \\  \Lambda \neq 1^n\ea}} \pm  \frac{k!}{\prod_{i=1}^k m_i(\Lambda) !} \dim A_{\Lambda} 
\]

Since $\dim A_{2i+1}=1$,
the dimension of a subspace $A_\Lambda =
A_{\Lambda_1}\ast \ldots \ast A_{\Lambda_k}\subset \Q[ \S_n]$ is
\[
\dim A_\Lambda = \frac{n!}{\prod_{i=1}^k \Lambda_i!}
\qquad |\Lambda|=n \ .
\]
Any subspace $A_{\Lambda_i}=A_{2\lambda_i+1}$ is in the intersection of $\lambda_i$ elementary subspaces
$A_{1^i31^j}$ 
hence the subspace $A_\Lambda$
lays in the intersection of $\sum_{i=1}^k \lambda_i$ subspaces.
The summation over partitions  $\Lambda$ having only odd parts is equivalent to the summation over partitions $\lambda$ (skipping
the empty diagram $\lambda=0$ corresponding to $\Lambda\neq 1^n$). The alternating sign factor is 
fixed by the number of intersections
$\sum_{i=1}^k \lambda_i=(|\Lambda|-k)/ 2$. 
These observations  together yields 
the final expression eq.(\ref{my}).

\section{Counting snakes and pre-plactic algebra}

The term ``snake" for an alternating permutation was coined by
 Vladimir Arnold \cite{Arnold}.

\begin{deff} An alternating permutation or {\it snake} is a permutation
  $$\sigma=\left(\ba{cccc}1 &2 &\ldots& n \\
	\sigma_1&\sigma_2 & \ldots& \sigma_n\ea\right)\in \S_n \qquad \mbox{such that}\qquad \sigma_1>\sigma_2<\sigma_3>\ldots \ .$$
		The word of the permutation $\sigma$ is the word
	$\sigma_1 \sigma_2 \ldots\sigma_n$.
	We denote by $E_n$ the number of snakes,
	$E_n=\#\{\sigma\in \S_n| \sigma_1>\sigma_2<\sigma_3>\ldots\}$. 
\end{deff}

The alternating permutations are counted by D\'esir\'e Andr\'e \cite{Andre}. Their number is generated
by the series $$y(x):=\sum_{n\geq 0} E_n\frac{x^n}{n!}=\sec x + \tan x$$
satisfying the equation\footnote{For a survey on alternating permutations see \cite{Stanley}} 
\beq
\label{sny}
 2y'=y^2 +1 \qquad y(0)=1 \ .
\eeq

 Arnold also  referred to  the numbers $E_n$ as to   Bernoulli-Euler numbers for the following reason:
the Euler numbers $E_{2n}$ appear from the expansion of
$$\sec x= \sum_{n\geq 0} E_{2n}\frac{x^{2n}}{(2n)!}$$
whereas the {\it tangent numbers} $E_{2n+1}$
$$
\tan x = \sum_{n\geq 0} E_{2n+1}\frac{x^{2n+1}}{(2n+1)!}
$$
are related to the Bernoulli numbers $E_{2n+1}=B_{2n}\frac{4^n(4^n-1)}{2n}$. Arnold introduced ``snakes" for other Coxeter groups
than $A_n$, their numbers 
are topological invariants of bifurcation diagrams \cite{Arnold}.

We found an expression of the snakes numbers $E_n$
through summation over Young partitions which to the best of
our knowledge is new.

\begin{prop}
The number of snakes 
 is given by the formula
\beq
\label{BE}
\frac{E_{n+1}}{n!}=\sum_{k > 0}
 \sum_{\Lambda \vdash n }  \frac{(-1)^{(|\Lambda|-k)/2}k!}{\prod_{i=1}^k  \Lambda_i ! \prod_{i=1}^k{m_i(\Lambda) !}}\qquad \mbox{with} \qquad \Lambda_i= 2 \lambda_i + 1
\eeq
where the sum is over all partitions $\Lambda$ of $n$ with 
only odd parts.
\end{prop}

{\bf Proof.} Differentiating the generating series $y(x)$ one gets
\[
y^\prime(x)=\sum_{n=0}^\infty E_{n+1} \frac{x^n}{n!}
=\frac{1}{1-\sin x}
= \sum_{k\geq 0} \sin^k x=
\sum_{k\geq 0} \left(\sum_{\lambda=0}^\infty 
(-1)^{\lambda} \frac{x^{2\lambda +1}}{(2\lambda+1)!} \right)^k
\]
The result follows by the binomial expansion and
equating the coefficients of $x^n$.$\Box$

We are now ready to prove the main result of this paper
namely the snakes numbers are the coefficients in
 the generating series for the $\S$-module $\PP$.

\begin{thm}
The generating series of the $\S$-module $\PP=\bigoplus_{n\geq 0} \PP(n)$ 
is given by
\[
\sum_{n\geq 0}\dim \PP(n)\frac{x^n}{n!} 
= 
\sum_{n=0}^\infty E_{n+1}\frac{x^n}{n!} \ .
\]
\end{thm}
{\bf Proof.} The homogeneous degree $\PP(n)$ of the pre-plactic algebra is 
 $\PP(n)=\Q[\S_n]/\mathfrak I(n)$ hence 
\beq
\label{match}
\dim \PP(n) = n!- \dim \mathfrak I(n) \ .
\eeq
The dimension of the ideal  $\dim \mathfrak I(n)$ was
found in  eq. (\ref{my}) from where we get 
\[
\dim \PP(n)= 
{\sum_{ k >0 }}\,\,
 \sum_{{\Lambda}\vdash
 n } (-1)^{(|{\Lambda}|-k)/2}
\frac{k!}{ \prod_{i=1}^k{m_i({\Lambda}) !}}
\frac{n!}{\prod_{i=1}^k  {\Lambda}_i!}  \ .
\]   
Note that by removing the restriction $\Lambda\neq 1^n$ in eq. (\ref{my}) we incorporated the extra  term $n!$ in $\dim \PP(n)$.
We now realize that the expression for $\dim \PP(n)$
  coincides with the formula eq. (\ref{BE}) for the snakes number 
$E_{n+1}$  
which readily implies the result $\dim \PP(n)=E_{n+1}$.

\section{Hopf structure and Hilbert-Poincar\'e Series}

{\bf Tensor algebra.}
The free associative algebra $T(V)$ generated by $D$-dimensional
vector space $V$ is naturally a  $GL(V)$-module.
Its character reads
\[
ch_{T(V)}(x)=\frac{1}{1-(x_1+ \ldots + x_D)}= \sum_{n\geq 0} (x_1+ \ldots + x_D)^n = \sum_\lambda s_\lambda(x) f^\lambda \ .
\]

By Schur-Weyl duality the image of 
$T(V)=\bigoplus_{n\geq 0} V^{\otimes n}$ is the $\S$-module
$\S=\bigoplus_{n\geq 0} \S_n$
with generating series
\beq
\label{cs}
g_{\S}(x)=\frac{1}{1-x} = \sum_{n\geq 0} n!\frac{x^n}{n!}=
\exp \ln \frac{1}{(1-x)} = \exp \sum_{k\geq 1} \frac{x^k}{k}=
\exp \sum_{k\geq 1} (k-1)!\frac{x^k}{k!}
\eeq
where $(k-1)!$ is the dimension of the 
space $Lie(k)\subset \Q[\S_k]$ of primitive Lie elements.

\noindent
{\bf Parastatistics algebra $PS(V)$.}
The universal enveloping algebra (UEA) of the  two step nilpotent Lie algebra 
$PS(V):=U(V\oplus\wedge^2 V)$ arises in the parastatistics Fock spaces \cite{D-VP} and it will be referred to as the {\it parastatistics algebra} $PS(V)$. The algebra
\[
PS(V) = T(V)/([[V,V],V]) 
\]
is a factor algebra of the tensor algebra $T(V)$ by the two-sided
ideal generated in $Lie(V)_3:=[[V,V],V]$ which is also a Hopf ideal
for the standard Hopf structure on $T(V)$.
The character of the $GL(V)$-module $PS(V)$ is also representable
as the sum of all Schur polynomials
\beq
\label{sch}
ch_{PS(V)}(x)=\prod_{i=1}^D \frac{1}{(1-x_i)}
\prod_{1\leq i<j\leq D} \frac{1}{(1-x_i x_j)} = \sum_{\lambda} s_\lambda(x_1, \ldots, x_D) \ .
\eeq
The latter identity implies 
that in the decomposition of $PS(V)$ into irreducible representation
every Schur module  apprears once and  exactly once,
$$PS(V)= \bigoplus_{n\geq 0} PS(n)\otimes_{\S_n} V^{\otimes n}
\cong \bigoplus_\lambda S_\lambda(V) \ .$$

The $\S$-module $PS=\bigoplus_{n\geq 0} PS(n) \cong \bigoplus_\lambda S_\lambda$ is the Schur-Weyl dual of the  $GL(V)$-module 
$PS(V) $. The dimension  of the irreducible $\S_n$-module $S_\lambda$
is the number $f_\lambda$ of  the Standard Young tableaux of shape
$\lambda$, $f_\lambda:=\dim S_\lambda$.
Hence the generating series of the Standard Young Tableaux
yields the generating series of the $\S$-module $PS$
\beq
\label{css}
g_{PS}(x)=\exp \left({x+ \frac{x^2}{2}}\right)=\sum_{n\geq 0} x^n\sum_{\lambda \vdash n} f_\lambda  \ .
\eeq

\noindent
{\bf Plactic algebra $\Pl(V)$.}
The plactic algebra $\Pl(V)$ turns out to be a ``crystal" limit
of a quantum deformation $PS_q(V)$ on the parastatistics algebra
$PS(V)$ for the singular value $q=0$ of the deformation parameter,
$\Pl(V)= PS_0(V)$ \cite{LP2}. The deformed algebra $PS_q(V)$
is a module of the quantum UEA $U_q \mathfrak gl_D$.
 The irreducible $U_q \mathfrak gl_D$-modules are the quantum Schur modules $S_\lambda^q(V)$ having the same dimensions as the Schur modules $S_\lambda(V)$. Hence the character formula (\ref{sch})
gives also the multi-graded dimensions of the deformed algebra $PS_q(V)$  and its specialization: the plactic algebra $\Pl(V)$
\[
ch_{PS(V)}(x)=ch_{PS_q(V)}(x)=ch_{\Pl(V)}(x)
\]

By quantum Schur-Weyl duality to the  $U_q \mathfrak gl_D$-module 
 $$PS_q(V)=
\bigoplus_{r\geq 0} PS_q(r)\otimes_{\hr} V^{\otimes r}
\cong \bigoplus_\lambda S^q_\lambda(V)  $$ one attaches the Hecke algebra modules 
$PS_q\cong \bigoplus_\lambda S_\lambda^q$.

The Poirier-Reutenauer algebra can be seen as the $q=0$ limit  
of direct sum of Hecke algebra modules \cite{LP3}.
 The plactic classes (\ref{plcl}) close a commutatuve
Hopf subalgebra (with respect to the shuffle) of the Malvenuto-Reutenauer algebra . Thus the generating series (\ref{css}) is also counting the plactic classes which span the Poirier-Reutenauer algebra $Tab^\prime$
\[
g_{PS}(x)= g_{PS_q}(x)=g_{Tab}(x) \ . 
\]

The  series eq.(\ref{css}) is a truncation of the generating series  
$\exp \sum_{k\geq 1} \frac{x^k}{k}$  of $\S$ 
reflecting  the factorization of the free Lie algebra
to the two-step nilpotent Lie algebra.
As a digression we now prove an identity which is one more truncation of the generating series  of $\S$
\begin{lemma}
 The snake numbers satisfy the  identity
\beq
\label{ide}
\sum_{n\geq 0} E_{n+1}\frac{x^n}{n!} =
\exp {\sum_{n\geq 1} E_{n-1}\frac{x^n}{n!}}  \ .
\eeq
\end{lemma}
{\bf Proof.}  The number of alternating permutations
$E_n$ is generated by the series $y(x)=\sec x + \tan x$
which is a solution of the equation $2y^\prime = y^2+1$. The identity is equivalent
to the following integro-differential equation
\beq
\label{intdif}
 y^{\prime}(x)  =\sum_{n\geq 0} E_{n+1} \frac{x^n}{n!}=
\exp {\int_0^{x} y(t)dt}
\eeq
for the generating series $y(x)=\sum_{n\geq 0} E_n\frac{x^n}{n!} =\sec x + \tan x$ of  the snakes numbers $E_n$.

Indeed the differention of the eq. (\ref{intdif}) yields
\[
 y^{\prime\prime}(x) = y(x)\exp {\int_0^{x} y(t)dt}\qquad
\Rightarrow \qquad
y^{\prime\prime} =yy^\prime  \qquad \Leftarrow \qquad (2 y^\prime)^{\prime} = (y^2+1)^\prime
\]
which is the same as the differentiation of the equation 
(\ref{sny}) for the generating function $y(x)$ of the snakes numbers.
$\Box$

\noindent
{\bf Pre-plactic algebra $\PP$ and  pseudo-plactic algebra $\PP(W)$.}
The left hand side of the   identity (\ref{ide}) 
is 
the generating series of the pre-placitic  algebra $\PP$.
Then the natural interpretation of the numbers in the exponent 
of the identity (\ref{ide}) is as dimensions of 
the space of primitive Lie  elements  ${\prPP}(n)$ 
of the Hopf algebra $\PP$
\[
g_{\PP}(x)=\sum_{n\geq 0} \dim \PP(n)\frac{x^n}{n!}
=\exp \left( {\sum_{n\geq 1} \dim \prPP(n)\frac{x^n}{n!}}\right) \ .
\]
The Malvenuto-Reutenauer algebra $\S$ is a Hopf algebra which
is not cocommutative and it is not an UEA of its Lie algebra.
The primitive elements of $\S$ with respect to $\Delta$ are indexed by the connected permutations \cite{DHT, PR}.
However among the primitives of $\S_n$ we single out the Lie
primitives which are Lie elements $Lie(n)$, $\dim Lie(n)=(n-1)!$.
Some of these primitive elements are in the ideal of $\PP(n)$
but the identity (\ref{ide}) that the dimension 
 of the non-vanishing
Lie elements  ${\prPP}(n)$ is given  by the snake number $E_{n-1}$.
To set the  bijection between alternating permutations and some 
commutative subalgebra of (pre-plactic) classes in $(\S,\ast^\prime, \Delta^\prime)$ is an interesting open problem that we hope to address in a near future.

{In summary, snakes numbers yields the dimensions of both
the factorizable and nonfactorizable (primitive) elements
of $\PP$.  Hence the Hopf algebra $\PP$ is a categorification
of the snakes numbers in the same fashion as 
the Malvenuto-Reutenauer Hopf algebra is a categorification
of the factorial numbers $n!$. 



\begin{thm} The character of the $GL(V)$-module $\PP(W)$
is given by
\beq
ch_{\PP(W)}(x_1, \ldots, x_D) =
{\prod_{n=1}^{\infty} \prod_{1\leq i_1<\ldots < i_n\leq D}
{(1-x_{i_1} x_{i_2}\ldots x_{i_n})^{-E_{n-1}}} } \ .
\eeq
It coincides with the character of 
the $U_q \mathfrak{gl}_D$-module $\PP_q(W)$
\[ 
ch_{\PP(W)}(x_1, \ldots, x_D) =ch_{\PP_q(W)}(x_1, \ldots, x_D)  \ .
\]
\end{thm}
{\bf Proof.}
The primitive Lie elements $[[a,b],a]$ and $[b,[a,b]]$, $1\leq a<b\leq D$ with repeating letters are  generators of the ideal $\I(W)$ of $\PP(W)$. Hence any bracketing of Lyndon words with repeating letter is in the ideal $\I(W)$. Therefore the non-vanishing  primitive Lie elements in $\prPP(V)_n$ stem from Lyndon words
of length $n$  without repeating letter (from the alphabet 
$\{x_1, \ldots, x_D\}$).
Then the Poincar\'e-Birkhoff-Witt theorem implies the expression for the character $ch_{\PP(W)}(x_1, \ldots, x_D)$. The
character $ch_{\PP_q(W)}(x_1, \ldots, x_D)$ is the same
since $\PP_q(W)$ is a deformation of $\PP(W)$.

\begin{cor}
The Hilbert series of the pseudo-plactic algebra $\PP(W)$
and $\PP_q(W)$ reads
\beq
H_{\PP(W)}(t)= 
\prod_{n=1}^{\infty} {(1-t^n)^{-E_{n-1}{\tiny \left(\ba{c}\dim V\\n \ea\right)}} } \ . 
\eeq
\end{cor}

\subsection*{Outlook and perspectives}
We designed the quantum pre-plactic algebra $\PP_q$ as a tool for proving the conjecture of Krob and Thibon about the isomorphism between the quantum pseudo-plactic algebra and the diagonal subalgebra of $\C[GL_q(V)]$ \cite{KT}. The sketch of this approach is given in the proceedings \cite{P}.

The pre-plactic algebra $\PP$ turn out to be an object with a rich combinatorial structure. An interesting open problem is to construct bijection between  snakes and the  primitive Lie elements in $\PP$. Another problem is to clarify the connection between the 
pre-plactic algebra $\PP$ and the non-commutative trigonometric
functions in $\bf FQSym$ \cite{NCsnakes} realization of snakes
of type $A$. More distant goal will be to study the diagonal subalgebra of functions on quantum groups
$\C[SO_q(n)]$ and $\C[Sp_q(2n)]$  and how they are related to the Arnold's snakes of type $B$ and $D$.

\subsection*{Acknowledgement}
It is my pleasure to thank Anton Alekseev,
Peter Dalakov, Tekin Dereli, Michel Dubois-Violette,
 G\'erard Duchamp, Ludmil Hadjiivanov 
and Oleg Ogievetsky for their interest in that work and for many inspiring discussions. 
This work 
has been supported  by Grant DFNI T02/6 of the Bulgarian National Science Foundation, TUBITAK 2221 program
and  a Ko\c c University grant.


\begin{thebibliography}{15}

\bibitem{AS}
Aguiar, Marcelo, and Frank Sottile. "Structure of the Malvenuto–Reutenauer Hopf algebra of permutations." Advances in Mathematics 191.2 (2005): 225-275.

\bibitem{Andre}
D. Andr\'e, D\'eveloppement de sec x  et  de tang x, C. R. Math. Acad. Sci. Paris 88 (1879), 965--979.
\bibitem{Arnold}
V. I. Arnol'd, The calculus of snakes and the combinatorics
of Bernoulli, Euler and Springer numbers of Coxeter groups,
Uspekhi Mat. Nauk, 1992, Volume 47, Issue 1(283), 3--45.

\bibitem{D-VP} M. Dubois-Violette and T. Popov,
         Homogeneous algebras, statistics and combinatorics.
          {\it Lett. Math. Phys. \/ \bf 61} (2002), 159-170. 


\bibitem{DHT}Duchamp, G\'erard, Florent Hivert, and Jean-Yves Thibon. "Noncommutative symmetric functions VI: free quasi-symmetric functions and related algebras." International Journal of Algebra and computation 12.05 (2002): 671-717.

\bibitem{NCsnakes} Josuat-Verges, Matthieu, Jean-Christophe Novelli, and Jean-Yves Thibon. "The algebraic combinatorics of snakes." Journal of Combinatorial Theory, Series A 119.8 (2012): 1613-1638.

\bibitem{K}D.E.Knuth, Permutation Matrices and Generalized Young Tableaux {\em Pacific Journal \/\bf 34}(1970),
709-727.    

\bibitem{KT}Krob, Daniel, and Jean-Yves Thibon. "Noncommutative symmetric functions IV: Quantum linear groups and Hecke algebras at q= 0." Journal of Algebraic Combinatorics 6.4 (1997): 339-376.



\bibitem{LP2} {J.-L. Loday and T. Popov},
 Parastatistics Algebra, Young Tableaux and Super Plactic Monoid. 
 {\em	International Journal of Geometric Methods in Modern Physics\/ \bf 5} (2008),  1295-1314.

\bibitem{LP3} J.-L. Loday, T. Popov.
Hopf Structures on Standard Young Tableaux.
{\em Proceedings 
 "Lie Theory and Its Applicatioins"},
  Vl. Dobrev ed.,
  AIP conference series {Vol. \bf 1243}(2010), 265--275.

\bibitem{MR} C. Malvenuto and C. Reutenauer, Duality 
between Quasi-Symmetric Functions and the Solomon Descent Algebra.
{ \em Journal of Algebra {\bf 177}} (1995), 967--982.        
\bibitem{PR} S.  Poirier and C. Reutenauer,  Alg\`ebres 
de Hopf de tableaux. {\em Ann. Sci. Math. Qu\' ebec \/ \bf 19 } (1995), 
79-90.

\bibitem{oleg}
O.
Ogievetsky, Uses of quantum spaces. {\em Contemporary Math. {\bf 294} }
(2002), 161-232.

\bibitem{P} Popov, Todor. "Quantum Plactic and Pseudo-Plactic Algebras." International Workshop on Lie Theory and Its Applications in Physics. Springer Singapore, 2015.

\bibitem{Stanley}
Stanley, Richard P. "A survey of alternating permutations." Contemp. Math 531 (2010): 165-196.

\bibitem{S}
Solomon, Louis. "A Mackey formula in the group ring of a Coxeter group." Journal of Algebra 41.2 (1976): 255-264.
\end{thebibliography}
\end{document}